\documentclass[11pt]{article}
\usepackage{amssymb}
\usepackage{graphicx}
\usepackage{latexsym}
\def\part#1{\frac{\partial\phantom{q}}{\partial#1}}

\newenvironment{rmk}{\begin{trivlist}\item[]{\bf Remark:} }
{\end{trivlist}}

\newenvironment{ex}{\begin{trivlist}\item[]{\bf Example:} }
{\end{trivlist}}
\newenvironment{exs}{\begin{trivlist}\item[]{\bf Examples:} }
{\end{trivlist}}
\newenvironment{prf}{\begin{trivlist}\item[]{\bf Proof:} }
{\hfill $\Box$ \end{trivlist}}
\newenvironment{lemprf}{\begin{trivlist}\item[]{\bf Proof:} }
 {\end{trivlist}}
\newtheorem{thm}{Theorem}
\newtheorem{definition}{Definition}
\newtheorem{prp}[thm]{Proposition}
\newtheorem{lemma}[thm]{Lemma}

\def\End{\mathop{\rm End}\nolimits}
\def\Ext{\mathop{\rm Ext}\nolimits}

\def\deg{\mathop{\rm deg}\nolimits}

\def\rk{\mathop{\rm rk}\nolimits}

\def\td{\mathop{\rm td}\nolimits}

\def\ch{\mathrm {ch}}

\newcommand{\R}{\mathbf{R}}
\newcommand{\C}{\mathbf{C}}

\newcommand{\Z}{\mathbf{Z}}

\newcommand{\CP}{{\mathbf P}}
\newcommand{\PP}{{\mathbf {\rm P}}}

\newcommand{\barT}{{\overline T}}
\newcommand{\barD}{{\overline D}}
\newcommand{\barE}{{\overline E}}
\begin{document}
\title{Generalized holomorphic bundles\\
 and the B-field action}
 \author{Nigel Hitchin\\[5pt]}
 \maketitle
\centerline{{\it Subject classification}: {Primary 53D18}}
\section{Introduction}
On a manifold with a generalized complex structure one can define the notion of a generalized holomorphic bundle. In this paper we study the simplest case, where the generalized complex structure is just an ordinary complex structure. Here, as noted by Gualtieri \cite{MG1},  the definition of a generalized holomorphic bundle yields a {\it co-Higgs bundle}. This is an adaptation of Simpson's Higgs bundles, but with the tangent bundle replacing the cotangent bundle. It consists of a holomorphic vector bundle $V$ and a holomorphic section $\phi$ of $\End V\otimes T$ satisfying $\phi\wedge\phi=0$.

The theory of Higgs bundles is now well-developed, especially for curves of genus greater than one and so, to point out the similarities and differences,   we present initially a few features and constructions of co-Higgs bundles on a complex manifold $M$ of arbitrary dimension, in particular when $M$ is a projective space. Our main point, however, is to investigate the new feature that generalized geometry brings to bear on the subject, which is the B-field symmetry. A real closed $(1,1)$-form $B$ can be considered as a symmetry of a complex structure regarded as a generalized complex structure, and as a consequence it transforms one co-Higgs bundle to another. We observe that stability is preserved, and that if the Dolbeault cohomology class is trivial, the action is a complex gauge transformation. In the simplest case of bundles on the projective line we show that this action is essentially integrating Nahm's equations. Reinterpreting a co-Higgs bundle \`a la Simpson as a sheaf on the total space $TM$ of the tangent bundle, the B-field action  is shown to be tensoring the sheaf by a line bundle over $TM$ defined by $B$.

We also consider  a twisted version of this defined by a 1-cocycle of closed $(1,1)$-forms. In generalized geometry this consists of replacing the bundle $T\oplus T^*$ by a more general Courant algebroid. Here we find that the spectral data consists of a sheaf over a holomorphic gerbe over $TM$ defined by the twisting cocycle. 

From its definition, a generalized holomorphic bundle  defines an elliptic complex analogous to the Dolbeault complex, which we call here the $\barD$-complex. Since a sheaf over a gerbe has no obvious cohomology groups, we investigate how to compute the $\barD$-cohomology in the generic case. This involves the  sheaf restricted to the zero section of $TM$, where the gerbe is canonically trivial.

\section{Generalized complex geometry}\label{gen}
Although we shall mainly be discussing ordinary complex manifolds, it is the viewpoint offered by the context of generalized geometry which drives this work. This setting consists of replacing the tangent bundle $T$ of a manifold  by $T\oplus T^*$, equipped with the natural indefinite inner product $(X+\xi,X+\xi)=i_X\xi$ and replacing the Lie bracket of vector fields with the skew-symmetric Courant bracket on sections of $T\oplus T^*$:
$$[X+\xi,Y+\eta]=[X,Y]+{\mathcal L}_X\eta-{\mathcal L}_Y\xi-\frac{1}{2}d(i_X\eta-i_Y\xi).$$
A closed 2-form $B$ then defines an action
$X+\xi\mapsto X+\xi+i_XB$ which preserves the inner product and commutes with the Courant bracket. This we call the B-field action. 
 
 A {\it generalized complex structure} \cite{NJH1} on a manifold of dimension $2m$ is defined to be an isotropic rank $2m$ subbundle $E^{1,0}\subset (T\oplus T^*)\otimes \C$ with $E^{1,0}\cap \barE^{1,0}=0$ and such that sections of $E^{1,0}$ are closed under the Courant bracket. We put $E^{0,1}= \barE^{1,0}$. The subbundles $E^{1,0}, E^{0,1}$ are the $\pm i$ eigenspaces of a complex structure $J$ on the real vector bundle $T\oplus T^*$.
 
 A complex manifold with local holomorphic coordinates $z_1,\dots,z_m$ defines a generalized structure where $E^{1,0}$ is spanned by 
  $$\frac{\partial}{\partial \bar z_1}, \frac{\partial}{\partial \bar z_2},\dots , dz_1, dz_2,\dots$$
 If $B$ is now a real closed form of type $(1,1)$, then the interior product with $\partial/\partial \bar z_i$ is the $(1,0)$-form $B_{\bar i j}dz_j$ so that the B-field action preserves $E^{1,0}$. In the context of generalized geometry this is a symmetry of the structure. In fact, given a real function $f$, $B=dd^c\!f$ defines a symmetry analogous to a Hamiltonian flow on a symplectic manifold.
 
 We shall need a twisted version of this: given a good open covering $\{U_{\alpha}\}$ and a \v Cech 1-cocycle $B_{\alpha\beta}$ of closed 2-forms we can identify $T\oplus T^*$ over $U_{\alpha}$ with $T\oplus T^*$ over $U_{\beta}$ by the action of $B_{\alpha\beta}$. The inner product and Courant bracket are preserved so this gives us a bundle $E$, locally modelled on $T\oplus T^*$, with both these features. Since  the B-field acts trivially on $T^*\subset T\oplus T^*$ we have an extension (an ``exact Courant algebroid")
 $$0\rightarrow T^*\rightarrow E\rightarrow T\rightarrow 0.$$
 The definition of generalized complex structure extends directly. In particular,  if $B_{\alpha\beta}$ is a cocycle of real closed $(1,1)$ forms on a complex manifold, then it defines a bundle $E$ with a generalized complex structure. 
 
 In \cite{MG1} Gualtieri introduced the concept of a {\it generalized holomorphic bundle} on a generalized complex manifold. Note first that, for any function $f$,  $df$ is a section of $T^*\subset E$ and with $E\otimes \C=E^{1,0}\oplus E^{0,1}$, we define $\bar\partial \!f$ to be the component in $E^{0,1}$. It is easy to see that  for an ordinary complex structure this is the usual $\bar\partial$-operator. 
  
 \begin{definition} \label{Dbar} A generalized holomorphic bundle on a generalized complex manifold  is a vector bundle $V$ with a differential operator $\barD:C^{\infty}(V)\rightarrow C^{\infty}(V\otimes E^{0,1})$ such that for a smooth function $f$ and section $s$ 
\begin{itemize}
\item
${\overline D}(fs)=\bar\partial \!f s+f\barD s$
\item
${\overline D}^2: C^{\infty}(V)\rightarrow C^{\infty}(V\otimes \Lambda^2E^{0,1})$ vanishes.
\end{itemize}
\end{definition}
The definition is a straight analogue of Malgrange's result that a $\bar\partial$-operator $\bar\partial_A$ on a vector bundle satisfying  $\bar\partial_A^2=0$ defines a holomorphic structure. In our case the second condition involves an extension of the operator $\barD$ to higher degree exterior powers of $E^{0,1}$. The simplest way to describe this is to note that the inner product identifies $E^{0,1}$ with the dual of the maximal isotropic subspace $E^{1,0}$ so we can describe an element of $\Lambda^p E^{0,1}$ as an alternating multilinear form on $E^{1,0}$. One then defines a generalized $\bar\partial$-operator on $\Lambda^p E^{0,1}$ just as for the exterior derivative, replacing the Lie bracket by the Courant bracket. Although the Courant bracket does not satisfy the Jacobi identity in general, it does on sections of  $E^{1,0}$ because of the  isotropy condition. Given the generalized $\bar\partial$-operator, it is straightforward to extend $\barD$.
 
 These issues hardly concern us, however, since we will be looking just at an ordinary complex structure. Then   $E^{0,1}= {\overline T}^*\oplus T$ (where now $T$ is the holomorphic tangent bundle) and a  generalized holomorphic bundle is defined by an operator (we use the summation convention) 
$$\barD s=\bar\partial_A s+\phi s=\left(\frac{\partial s}{\partial \bar z_{j}}+A_{\bar j}s\right) d\bar z_{j}+\phi^{k} s \frac{\partial }{\partial z_{k}}$$
and $\barD^2=0$ implies 
\begin{itemize}
\item
$\bar\partial_A^2=0 \in C^{\infty}(\End V\otimes \Lambda^2\barT^*)$
\item
$ \bar\partial_A \phi=0 \in C^{\infty}(\End V\otimes T\otimes \barT^*)$
\item 
$ \phi\wedge\phi=0\in C^{\infty}(\End V\otimes \Lambda^2T).$
\end{itemize}
The first condition, by Malgrange's result, gives $V$ the structure of a holomorphic vector bundle, the second says that $\phi$ is a holomorphic section of $\End V\otimes T$, and the third  is an algebraic condition on $\phi$. Since
$$\phi\wedge\phi=\frac{1}{2}[\phi^{j},\phi^{k}]\frac{\partial }{\partial  z_{j}}\wedge \frac{\partial }{\partial  z_{k}}$$
this ``integrability" condition is   $[\phi^{j},\phi^{k}]=0$. 
 
 With  $T$ replaced by $T^*$, this is precisely Simpson's definition of {\it Higgs bundle} \cite{S}. We shall call these structures --  generalized holomorphic bundles for an ordinary complex structure -- {\it co-Higgs bundles}.
 
 \section{Co-Higgs bundles}
 \subsection{Examples}
 Here are some simple examples of co-Higgs bundles:

\noindent 1. In rank one, a co-Higgs bundle is just a holomorphic line bundle $L$ together with a vector field $X$.

\noindent 2. If $V={\mathcal O}\oplus T$ then there is a canonical co-Higgs structure where $\phi(\lambda,X)=(X,0)$. 

\noindent 3. Take $V$ to be the trivial bundle of rank $k$ on $\CP^1$. The tangent bundle is the line bundle ${\mathcal O}(2)$ so $\phi$ is a $k\times k$ matrix with entries in $H^0(\CP^1,{\mathcal O}(2))$ or, more concretely, quadratic polynomials in an affine coordinate $z$. In one dimension the condition $\phi^2=0$ is trivially satisfied.

\noindent 4.  To get examples with more interesting Higgs fields we can look at projective spaces. Consider $\CP^m=\PP(W)$. The tangent bundle fits into the Euler sequence of holomorphic bundles 
$$0\rightarrow  {\mathcal O}\rightarrow         W\otimes {\mathcal O}(1)\rightarrow T\rightarrow 0$$
from which we obtain $W\cong H^0(\CP^m,T(-1)).$ Now let $X$ be a smooth $m$-dimensional projective variety and embed it in a high-dimensional projective space $\CP^N$ and take a generic projection $\pi$ to $\CP^m$, so that restricted to $X\subset \CP^N$ we have a finite surjective morphism. Then for any line bundle $L$ on $X$, the direct image sheaf $\pi_*(L)$ is a holomorphic vector bundle $V$ whose rank is the degree of $X$. Each section $s$ of ${\mathcal O}(1)$ on $\CP^N$, restricted to $X$, defines by multiplication a homomorphism 
$$H^0(\pi^{-1}(U),L)\rightarrow H^0(\pi^{-1}(U),L(1))$$
which defines a section $\psi_s\in H^0(\CP^m,\End V(1))$. Moreover, since  multiplication by $s$ and $t$ commute, the twisted endomorphisms $\psi_s,\psi_t$ commute. Sections pulled back from $\CP^m$ give multiples of the identity but in general, if $X$ lies in no hyperplane, we have an $(N-m)$-dimensional space of trace zero twisted endomorphisms with basis $\psi_1,\dots,\psi_{N-m}$.

Now take sections $Y_i$ of $T(-1)$ on $\CP^m$ and write
$$\phi=\sum_{i=1}^{N-m}\psi_i\otimes Y_i\in H^0(\CP^m,\End V\otimes T).$$
By construction, since the $\psi_i$ commute, we have $\phi\wedge\phi=0$.

The case where $N=3$ and $X$ is a quadric surface gives the Schwarzenberger bundles \cite{Sch} on $\CP^2$. A detailed study of these as co-Higgs bundles can be found in Steven Rayan's forthcoming Oxford thesis. 

\subsection{The spectral variety}
 Simpson in \cite{S1} reinterpreted a Higgs bundle, or more generally a Higgs sheaf, on $M$ in terms of a sheaf ${\mathcal S}$ on $\PP({\mathcal O}\oplus T^*)$ whose support is disjoint from the divisor at infinity.  This can be adapted immediately, replacing $T^*$ by $T$. This algebraic approach is essentially the idea that a local 1-form on $M$ acts on $V$ by contraction with the Higgs field $\phi\in H^0(M,\End V\otimes T)$. The condition $\phi^2=0$ means that this is a commutative action and so ${\mathcal O}(V)$ is a module, not only over local holomorphic functions on $M$, but also  over the symmetric algebra generated by $T^*$, or the functions on $TM$ polynomial in the fibre directions. The support of this sheaf is the spectral variety. 

In local terms we have 
$$\phi=\phi^{k}  \frac{\partial }{\partial z_{k}}.$$
The coordinates $z_1,\dots,z_m$ define in the standard way 
$$(z_1,\dots,z_m,y_1,\dots, y_m)\mapsto  y_{k}  \frac{\partial }{\partial z_{k}}$$
coordinates $y_1,\dots,y_m,z_1,\dots,z_m$ on $TM$. Then $y_k$ acts on ${\mathcal O}_M(V)$ as $\phi^k$, and this way we obtain a module over ${\mathcal O}_{TM}$. 
The {\it spectral variety} $S\subset TM$ is defined by the characteristic equations 
$\det (a_k(y_k-\phi^k(z)))=0$ for $a\in \C^m$ and is the support of the sheaf ${\mathcal S}$.  Each $a_k\phi^k$ then has a distinguished eigenvalue $a_ky_k$ on $S$ and, because the $\phi^k$ commute, there is a common eigenspace. 
When the spectral variety is smooth and the sheaf is a line bundle $L$ supported on $S$, the vector bundle is recovered as a direct image sheaf.: $V=\pi_*L'$ where $L'=K_SK_M^*L^*$.

We  consider now the spectral varieties for the examples in the previous section.

\begin{exs}

\noindent 1. A rank one co-Higgs bundle is  a  line bundle $L$ together with a vector field $X$.  Considering the vector field as a section $X:M\rightarrow TM$, then the spectral variety $S=X(M)\subset TM$.  
The sheaf is just the line bundle $L$ on $S\cong M$.

\noindent 2. When $V={\mathcal O}\oplus T$ and $\phi(\lambda,X)=(X,0)$ then $\phi^j\phi^k=0$ so  $y_jy_k=0$ in its action on ${\mathcal O}(V)$. In particular, the sheaf is zero away from the zero section of $TM$. The  sheaf  ${\mathcal S}$ is  the trivial line bundle on the  first order neighborhood of the zero section of $TM$, and this is the spectral variety.

\noindent 3. When $V$ is the trivial rank $k$ bundle on $\CP^1$, $\phi$ is a $k\times k$ matrix of sections of ${\mathcal O}(2)$ and the spectral variety is the  curve
$$\det (y-\phi(z))=0$$
in the total space of $\pi:{\mathcal O}(2)\rightarrow {\CP^1}$.  Generically this is a smooth curve  of genus $(k-1)^2$ and the sheaf is a line bundle of degree $k-k^2$. 

\noindent 4. In the final example we built our Higgs field on $\CP^m$ out of the direct  image of sections $s_1,\dots,s_{N-m}$ of $\pi^*{\mathcal O}(1)$ and sections $Y_i$ of $T(-1)$. Then
$$\sum_{i=1}^{N-m} s_i\pi^*Y_i$$
is a section of $\pi^*T$ on $X$ which maps $X$ to the spectral variety $S$ in $T\CP^m$. 

When $N-m=m+1$ we can take the $Y_i$ above to be a basis of sections of $T(-1)$ over $\CP^m$ and then there is a more geometrical description. If $X$ is an $m$-dimensional submanifold of $\CP^{2m+1}$ we remove a subspace $\CP^m\subset \CP^{2m+1}$ which misses $X$, then the complement $\CP^{2m+1}\setminus \CP^m$ is isomorphic to the total space of ${\mathcal O}^{m+1}(1)$ over $\CP^m$. Choosing a non-vanishing section gives a quotient bundle   ${\mathcal O}^{m+1}(1)/{\mathcal O}$ which, as in the Euler sequence, is the tangent bundle $T\CP^m$. Then $X\subset {\mathcal O}^{m+1}(1)$  maps in the quotient to the spectral variety $S$.

 In concrete terms, if $[z_0,\dots,z_m,w_0,\dots,w_m]$ are homogeneous coordinates on $\CP^{2m+1}$ then removing the subspace defined by $z_0=\dots = z_m=0$, we project to the point $[z_0,\dots,z_m]\in \CP^m$. The section $w_i=z_i$ for $ 0\le i\le m$  is non-vanishing and the map to $T\CP^m$ in affine coordinates $[1,z_1,\dots,z_m]$ is given by 
 $$\sum_{i=1}^m(w_i-z_iw_0)\frac{\partial}{\partial z_i}.$$

\noindent 5. A specific example of the previous case is  the Veronese embedding of $X=\CP^2$ in $\CP^5$, given in homogeneous coordinates by 
$[x_0,x_1,x_2]\mapsto [x_0^2,x_1^2,x_2^2,x_1x_2,x_2x_0, x_0 x_1]$. 
The projection  onto $\CP^2$ is
$[x_0,x_1,x_2]\mapsto [x_0^2,x_1^2,x_2^2]$ 
which is a  fourfold covering, branched over three lines $z_0=0,z_1=0,z_2=0$. The derivative of the map $X\rightarrow T\CP^2$ is injective except possibly where $x_0=0,x_1=0$ or $x_2=0$ and indeed the standard Veronese surface acquires singularities under this map, but if we apply an invertible  linear transformation of the form $(z,w)\mapsto (z,Aw)$ then the finite set of inequalities on the coefficients of the matrix $A$ to ensure smoothness can be generically satisfied, and we obtain an embedding of $\CP^2$ as a surface $S$ in $T\CP^2$. Choosing a line bundle $L$, this becomes  the spectral variety for a rank $4$ co-Higgs bundle on $\CP^2$. 
 \end{exs}
\subsection{Stability}

On a K\"ahler manifold $(M,\omega)$ one can define the degree $\deg L=c_1(L)[\omega]^{m-1}$ of a line bundle  $L$ 
 and the slope $\mu(V)=\deg \Lambda^kV/\rk(V)$ of a vector bundle $V$. The stability condition for Higgs bundles in \cite{S} is that the slope of any $\phi$-invariant torsion-free subsheaf should be less than the slope of $V$ and this adapts immediately to the co-Higgs situation. Stability provides the means to define moduli spaces rather than stacks and also links up with gauge-theoretic notions.

\begin{exs}

\noindent 1.  In dimension one, stability can only occur, for rank greater than one, in genus $g\le 1$, for contracting $\phi$ with a holomorphic differential gives a non-trivial endomorphism of $V$ commuting with $\phi$ which is impossible in the stable case.

\noindent 2. If $T$ is stable as a vector bundle and $\deg T>0$, then the canonical  co-Higgs bundle ${\mathcal O}\oplus T$  is stable.

 \noindent 3. If the spectral variety is reduced and irreducible then there are no invariant subsheaves and the co-Higgs bundle is automatically stable. 
 
 \end{exs}

The gauge-theoretic link is provided by Simpson's theorem  \cite{S0} that a stable Higgs bundle has a canonical Hermitian metric satisfying the Yang-Mills type equation
$$(F+[\phi,\phi^*])\wedge\omega^{m-1}=\lambda I \omega^m$$
where $\lambda$ is a constant (in fact $\lambda=-i\mu(V)$)
and $F$ is the curvature of the unique connection compatible with the holomorphic structure on $V$ and the Hermitian metric. In coordinates,
$$\omega^{j\bar k}(F_{j\bar k}+[\phi_j,\phi^*_{\bar k}])=\lambda.
$$ The methods readily adapt (these are moment map equations for the gauge group) to the co-Higgs bundle case where the equation is 
\begin{equation}
\omega^{j\bar k}F_{j\bar k}+\omega_{j\bar k}[\phi^j,(\phi^*)^{\bar k}]=\lambda
\label{eqns}
\end{equation}
We shall give another interpretation of this equation  next.  Define $\psi\in \Omega^{0,1}(\End V)$ by the interior product $\psi=i_{\phi}(i \omega)$, i.e.
$$\psi_{\bar k}=i\phi^j\omega_{j\bar k}.$$
Then since $\omega$ is closed and $\phi$ is holomorphic, $\bar\partial_A\psi=0\in \Omega^{0,2}(\End V)$. Moreover, since $\phi\wedge\phi=0$, it follows that $\psi\wedge\psi=0$. This means that $\bar\partial_{+}=\bar\partial_A+\psi$ satisfies $\bar\partial_{+}^2=0$ and defines a new holomorphic structure on $V$. Similarly replacing $\omega$ by $-\omega$ we get a holomorphic structure $\bar\partial_{-}$. Given a Hermitian metric on $V$, there are unique unitary connections $\nabla_+$ compatible with $\bar\partial_{+}$ and $\nabla_-$ compatible with $\bar\partial_{-}$.
\begin{prp} The Hermitian metric satisfies the equation (\ref{eqns}) if and only if 
$$\frac{1}{2}(F_+ +F_-)\wedge\omega^{m-1}=\lambda I \omega^m$$
where $F_+,F_-$ are the curvatures of the connections $\nabla_+,\nabla_-$.
\end{prp}
\begin{prf}
Let $\nabla$ be the Hermitian connection compatible with $\bar\partial_A$, with curvature $F$. Then $\nabla_+=\nabla+\psi-\psi^*$ and 
$$F_+=F+\nabla^{1,0}\psi-\nabla^{0,1}\psi^*-\psi\psi^*-\psi^*\psi.$$
Similarly
$$F_-=F-\nabla^{1,0}\psi+\nabla^{0,1}\psi^*-\psi\psi^*-\psi^*\psi.$$
Taking the sum gives the result when we expand
$$(\psi\psi^*)_{\bar k m}=(i\phi^j\omega_{j\bar k})(-i(\phi^*)^{\bar\ell}\omega_{\bar \ell m})=-\phi^j(\phi^*)^{\bar\ell}\omega_{j\bar k}\omega_{m \bar \ell}$$
and then 
$$\omega^{\bar k m}(\psi\psi^*)_{\bar k m}=-\phi^j(\phi^*)^{\bar\ell}\omega^{\bar k m}\omega_{j\bar k}\omega_{m \bar \ell}=-\phi^j(\phi^*)^{\bar\ell}\omega_{j \bar \ell}.$$
\end{prf}

\begin{rmk} These equations are suggestive of Yang-Mills-type equations in the case of a {\it generalized K\"ahler structure} \cite{MG2}. In this situation we have two commuting generalized complex structures $J_1=J$ and $J_2$ and the eigenspaces of $J_2$ split $E^{0,1}$ into two subbundles. These define a bihermitian metric -- two complex structures $I_+,I_-$ and a single Hermitian metric $g$ on $M$. 
A generalized holomorphic bundle with respect to $J$ then  defines two holomorphic structures on the bundle, one with respect to $I_+$ the other with respect to $I_-$. What we did above was the ordinary K\"ahler case, where the second generalized complex structure is the symplectic one defined by the K\"ahler form. 

One would then expect in general an equation
$$\frac{1}{2}(F_+\wedge \omega_+^{m-1} +F_-\wedge \omega_-^{m-1})=\lambda I {\mathrm{vol}}_g$$
to describe a stability condition on the generalized holomorphic bundle. 
\end{rmk}
\section{The B-field action}
\subsection{The action}
On a complex manifold we can pull back a holomorphic vector bundle by a holomorphic diffeomorphism to get a new one. But we saw that a closed real $(1,1)$-form $B$ can also be seen as a symmetry of a complex manifold when we consider it as a generalized complex manifold.

A co-Higgs bundle gave us a  generalized holomorphic bundle  defined by an operator
  $$\barD s=\left(\frac{\partial s}{\partial \bar z_{i}}+A_{\bar i}s\right) d\bar z_{i}+\phi^{j} s \frac{\partial }{\partial z_{j}}$$
 and the transform of this by $B$ is 
  $$\barD s=\left(\frac{\partial s}{\partial \bar z_{i}}+A_{\bar i}s+\phi^{j}B_{j\bar i}s\right) d\bar z_{i}+\phi^{j} s \frac{\partial }{\partial z_{j}}.$$

    More invariantly we write $i_{\phi}B\in \Omega^{0,1}(\End V)$ for the contraction of the Higgs field $\phi\in H^0(M,\End V\otimes T)$ with $B\in \Omega^{0,1}(T^*)$ and then we have a new $\bar\partial$-operator 
  $$\bar\partial_B=\bar\partial_A+i_{\phi}B$$
  on the $C^{\infty}$ bundle $V$. Just as in the previous section, where we used the {\it imaginary} closed $(1,1)$-form $\omega$, this defines a new holomorphic structure on $V$. Moreover, since $\phi$ commutes with $i_{\phi}B$, $\phi$ is still holomorphic with respect to $\bar\partial_B$ and we have a new co-Higgs bundle.
 
\begin{rmk} Note that if $U\subset V$ is a  subbundle holomorphic with respect to $\bar\partial_A$ and  also $\phi$-invariant, it is holomorphic with respect to $\bar\partial_B$. So stability is preserved by the B-field action.
\end{rmk}

 We consider now our examples in this light:

 \begin{exs} 
 
 \noindent 1. A rank one co-Higgs bundle is a line bundle $L$ together with a holomorphic vector field $X$.  Changing $\bar\partial_A$ -- the holomorphic structure of $L$ -- to $\bar\partial_B=\bar\partial_A+i_{X}B$ is just tensoring $L$ by the degree zero line bundle $L_B$ defined by the $\bar\partial$-operator $\bar\partial +i_{X}B$ on the trivial bundle.

 \noindent 2. Take the canonical Higgs bundle $V={\mathcal O}\oplus T$ and  $\phi(\lambda,X)=(X,0)$.
 Then in the direct sum decomposition 
 $$\bar\partial_A=\pmatrix {\bar\partial & 0\cr
 0 &\bar\partial}\qquad \phi=\pmatrix {0 & 1\cr
 0 &0}.$$
 Regarding $B$ as lying in $\Omega^{0,1}(T^*)$ the new holomorphic structure is 
  $$\bar\partial_B=\pmatrix {\bar\partial & B\cr
 0 &\bar\partial}$$
   which defines an  extension 
 $$0\rightarrow {\mathcal O}\rightarrow V\stackrel{\pi}\rightarrow T\rightarrow 0.$$
Note that there is still  a canonical Higgs field $v\mapsto \pi v\in T\subset V\otimes T$.
 \end{exs}

\begin{prp} If $B=\bar\partial \theta$, then the two holomorphic structures $\bar\partial_A$ and $\bar\partial_B$ are equivalent by a complex gauge transformation.
\end{prp}
\begin{prf}
 If  $B=\bar\partial \theta$ for  $\theta \in \Omega^{1,0}$  define $\varphi=i_{\phi}\theta$, a section of $\End V$ or  in coordinates  $\varphi=\phi^{i}\theta_{i}$. We observe  
$(\bar\partial_A\varphi)_{\bar i}=\phi^{j}B_{j\bar i}$ since $\phi$ is holomorphic and then
$$[\varphi,\bar\partial_A\varphi]_{\bar i}=[\phi^{k}\theta_{k},\phi^{j}B_{j\bar i}]=[\phi^{k},\phi^{j}]\theta_{k}B_{j\bar i}=0$$
since $\phi^j$ and $\phi^k$ commute. This means  (unusually for a non-abelian gauge theory) that if we exponentiate $\varphi$ to be an automorphism of the bundle $V$ then we have
 $$\exp (-\varphi)\bar\partial_A\exp \varphi=\bar\partial_A+i_{\phi}B.$$
\end{prf}

 \subsection{The spectral action}
Now consider the B-field action from the point of view of the spectral sheaf. Since $\phi$ is unchanged, the support of the sheaf ${\mathcal S}$ on $TM$ -- the spectral variety -- is unchanged. To understand the change in the sheaf, let $t$ denote the tautological section of $\pi^*T$ over the total space of $\pi:TM\rightarrow M$.  In local coordinates 
$$t= y_{k}  \frac{\partial }{\partial z_{k}}.$$
Pull back $B$ to $TM$ as a $(0,1)$-form with values in $\pi^*T$ and contract with $t$. Then $i_{t}B\in \Omega_{TM}^{0,1}$ defines a holomorphic structure on the trivial $C^{\infty}$ bundle -- a line bundle $L_B$.

\begin{prp} \label{sheafact} Let ${\mathcal S}$ be the spectral sheaf on $TM$ of a co-Higgs bundle on $M$. Then the B-field transform has  spectral sheaf ${\mathcal S}\otimes {\mathcal O}(L_B).$
\end{prp}
\begin{prf}
We pass from Dolbeault to \v Cech terminology. On a suitably small  open set $U_{\alpha}\subset M$ write $B=\bar\partial \theta_{\alpha}$. On $U_{\alpha}\cap U_{\beta}$, $\bar\partial (\theta_{\beta}-\theta_{\alpha})=0$ so $\theta_{\beta}-\theta_{\alpha}$ is a holomorphic one-form
$\sigma_{\alpha\beta}=h_idz_i$.  In \v Cech terms, the line bundle $L_B$ on $TM$ is defined by the transition functions $\exp(h_iy_i)$ on $U_{\alpha}\cap U_{\beta}$.

Now,  as in the previous section, $\varphi_{\alpha}=i_{\phi}\theta_{\alpha}$ exponentiates to a local automorphism which transforms $\bar\partial_B$ to the original $\bar\partial_A$. So over $U_{\alpha}$ the sheaf ${\mathcal S}_B$ is isomorphic to ${\mathcal S}$. The two automorphisms over $U_{\alpha}$ and $U_{\beta}$ are related over the intersection by 
$$\exp(\varphi_{\beta}-\varphi_{\alpha})=\exp (h_i\phi^i).$$ But $y_i$ acts on the sheaf over $TM$ by $\phi_i$ so this is tensoring the sheaf by ${\mathcal O}(L_B)$.
\end{prf}

\begin{exs}

\noindent 1. Consider the tautological Higgs bundle ${\mathcal O}\oplus T$.  The spectral sheaf  was the trivial bundle on the first order neighbourhood. If ${\mathcal O}_Z^{(1)}$ denotes the structure sheaf of  the first order neighbourhood of the zero section $Z$, there is an exact sequence of sheaves 
$$0\rightarrow {\mathcal O}_Z(N^*)\rightarrow {\mathcal O}_Z^{(1)}\rightarrow {\mathcal O}_Z\rightarrow 0$$
where $N$ is the conormal bundle.   In this case $N\cong T$ and  from the long exact cohomology sequence 
$$0\rightarrow H^1(Z,T^*)\rightarrow H^1(Z,{\mathcal O}^{(1)})\rightarrow H^1(Z,{\mathcal O})\rightarrow \cdots $$
we see that any two line bundles of degree zero on the neighbourhood differ up to isomorphism by an element of $H^1(Z,T^*)=H^1(M,T^*)$. These are realized by B-field actions.
\end{exs}

\noindent 2. More interesting is the one-dimensional case. Here, since $ H^1(\CP^1,T^*)=H^1(\CP^1,{\mathcal O}(-2))$ is one-dimensional, there is basically a unique B-field action. The corresponding line bundle $L_B$ on $T\CP^1$ appears in the twistor description of magnetic monopoles in $\R^3$ in \cite{HM1} -- it is the twistor transform of the solution to the $U(1)$ Bogomolny equations consisting of a constant Higgs field.  Take $V$ to be the trivial rank $k$ bundle on $\CP^1$ and, in an affine coordinate $z$, set 
$$\phi=(T_1+iT_2)+2iT_3 z+(T_1-iT_2)z^2$$
where the $T_i$ are $k\times k$ matrices.
We have a co-Higgs bundle and then, as described in \cite{HM2}, tensoring the spectral sheaf with the line bundle $L_B$  is, up to conjugation,  the evolution at time $t=1$ of the Nahm equations:
 $$\frac{dT_1}{dt}=[T_2,T_3],\quad \frac{dT_2}{dt}=[T_3,T_1],\quad \frac{dT_3}{dt}=[T_1,T_2]$$
 which are instrumental in solving the $SU(2)$ Bogomolny equations of charge $k$. This, as we have seen, is the B-field action. 
 
 \subsection{Twisted bundles}\label{twistb}
 Now  replace the generalized complex structure on $T\oplus T^*$ by a twisted version on the bundle $E$ defined by a 1-cocycle $B_{\alpha\beta}$ of closed real $(1,1)$-forms. What is a generalized holomorphic bundle now? The general definition (\ref{Dbar}) is the same -- a vector bundle $V$ with a differential operator $\barD$ -- but it demands a purely holomorphic description.  

If we consider $E$ as obtained by patching together copies of $T\oplus T^*$ then over each open set $U_{\alpha}$, $V$ has the structure of a co-Higgs bundle -- a holomorphic structure $A_{\alpha}$ and a Higgs field $\phi_\alpha$. On the intersection $U_{\alpha}\cap U_{\beta}$ these are related by the B-field action of $B_{\alpha\beta}$:

 \begin{equation}
 (A_{\beta})_{\bar i}=(A_{\alpha})_{\bar i}+\phi^{j}(B_{\alpha\beta})_{j\bar i}, \qquad (\phi_{\beta})^{j} = (\phi_{\alpha})^{j}.
 \label{af}
 \end{equation}

Consider first the case of $V=L$ a line bundle. Then, because $\End V$ is holomorphically trivial for all  of the local holomorphic structures, $\phi$ is a global holomorphic vector field $X$. So consider the $(0,1)$ form 
$$A_{\alpha\beta}=i_XB_{\alpha\beta}.$$
The $(1,1)$ form $B_{\alpha\beta}$ is closed  and $X$ is holomorphic so that $\bar\partial A_{\alpha\beta}=0$. Thus $\bar\partial +  A_{\alpha\beta}$ defines a holomorphic structure  $L_{\alpha\beta}$ on the trivial $C^{\infty}$ line bundle over the intersection $U_{\alpha}\cap U_{\beta}$. Since $B_{\alpha\beta}$ is a 1-cocycle we have $B_{\alpha\beta}=-B_{\beta\alpha}$ and $B_{\alpha\beta}+B_{\beta\gamma}+B_{\gamma\alpha}=0$ on threefold intersections. This means we have canonical isomorphisms 
$L_{\alpha\beta}\cong L^*_{\beta\alpha}$ and a canonical trivialization of $L_{\alpha\beta}L_{\beta\gamma}L_{\gamma\alpha}$ on $U_{\alpha}\cap U_{\beta}\cap U_{\gamma}$. As in \cite{NJH0}, this describes a holomorphic gerbe on $M$.  

We could also write $B_{\alpha\beta}=\bar\partial \theta_{\alpha\beta}$ and then $h_{\alpha\beta\gamma}= i_X(\theta_{\alpha\beta}+ \theta_{\beta\gamma}+ \theta_{\gamma\alpha})$ is holomorphic and $g_{\alpha\beta\gamma}=\exp h_{\alpha\beta\gamma}$ is a 2-cocycle which also defines the gerbe (see \cite{Bry}). 

The equation (\ref{af}) gives $\bar\partial$-closed forms $A_{\alpha}$ with 
$$A_{\beta}=A_{\alpha}+A_{\alpha\beta}$$
which means holomorphic line bundle structures $L_{\alpha}$ over $U_{\alpha}$ with isomorphisms 
$$L_{\beta}\cong L_{\alpha}L_{\alpha\beta}$$
compatible with the isomorphism in the definition of the gerbe. In the language of \cite{NJH0} this is a {\it trivialization} of the gerbe. 

\begin{ex} The simplest example is the Hopf surface $M=\C^2\setminus \{0\}/\Gamma$ where, with $\vert\lambda\vert\ne 1$,  $\Gamma$ is the infinite cyclic group generated by $(z_1,z_2)\mapsto (\lambda z_1, \lambda z_2)$. Since $H^2(M,{\mathcal O}^*)=\Z$, then from the 2-cocycle point of view, we see that  any holomorphic gerbe which is topologically trivial is holomorphically trivial, and so there is a generalized line bundle. 

As an example of the local line bundle description, cover $M$ by two open sets  $U_1$ where $z_1\ne 0$ and $U_2$ where $z_2\ne 0$ and take the real closed $(1,1)$-form 
$$B_{12}=\frac{dz_1}{z_1}\frac{d\bar z_2}{\bar z_2}-\frac{dz_2}{z_2}\frac{d\bar z_1}{\bar z_1}.$$
The tangent bundle is trivial and a vector field $X$ is of the form
$$X=a_1z_1\frac{\partial}{\partial z_1}+a_2z_2\frac{\partial}{\partial z_2}$$
and so
$$A_{12}=i_XB_{12}=a_1\frac{d\bar z_2}{\bar z_2}-a_2\frac{d\bar z_1}{\bar z_1}=A_2-A_1.$$
The $(0,1)$-form $A_1=a_2d\bar z_1/\bar z_1$ defines a holomorphic line bundle $L_1$ on $U_1$ and similarly  $A_2$ defines $L_2$ on $U_2$. These bundles are generally non-trivial -- the map $(z_1,z_2)\mapsto [z_1,z_2]\in \CP^1$ represents $M$ as an elliptic fibration and $L_1,L_2$ are non-trivial in general on the fibres. 
\end{ex}

With the 2-cocyle  approach to gerbes one  defines more generally a sheaf over the gerbe using a specific covering by 
\begin{definition}
A sheaf over a gerbe $\{g_{\alpha\beta\gamma}\}$ is defined by a pair $(\{{\mathcal S}_{\alpha}\},\{\varphi_{\alpha\beta}\})$ where ${\mathcal S}_{\alpha}$ is a sheaf of ${\mathcal O}$-modules on $U_{\alpha}$ 
and $\varphi_{\alpha\beta}:{\mathcal S}_{\beta}\rightarrow {\mathcal S}_{\alpha}$ is an isomorphism over $U_{\alpha}\cap U_{\beta}$ such that
\begin{itemize}
\item
$\varphi_{\alpha\alpha}$ is the identity on $U_{\alpha}$
\item
$\varphi_{\alpha\beta}=\varphi_{\beta\alpha}^{-1}$ on $U_{\alpha}\cap U_{\beta}$
\item
 $\varphi_{\alpha\beta}\varphi_{\beta\gamma}\varphi_{\gamma\alpha}$ is multiplication by $g_{\alpha\beta\gamma}$ on $U_{\alpha}\cap U_{\beta}\cap U_{\gamma}$.
 \end{itemize}
 \end{definition}

In the rank one case we could therefore replace ``trivialization of a gerbe" by  ``a line bundle over the gerbe", for  writing $A_{\alpha\beta}=\bar\partial f_{\alpha\beta}$, we can define
$$\varphi_{\alpha\beta}=\exp f_{\alpha\beta}$$
and this gives an isomorphism from $L_{\alpha}$ to $L_{\beta}$ over $U_{\alpha}\cap U_{\beta}$. One easily checks that this satisfies the conditions of the definition.

More generally consider the projection $\pi:TM\rightarrow M$ and the tautological section $t$ of $\pi^*T$. Given the cocycle $B_{\alpha\beta}$ on $M$ we can form 
$$A_{\alpha\beta}=i_{t}\pi^*B_{\alpha\beta}$$
and this defines a holomorphic gerbe ${\mathcal G}$ on $TM$. 

\begin{rmk} Note that on the zero section $Z\subset TM$, $t=0$ and hence $A_{\alpha\beta}=0$ so the gerbe is canonically trivial on $Z$.
\end{rmk}

Now suppose we have a generalized holomorphic bundle relative to the twisted generalized complex structure defined by the cocycle $B_{\alpha\beta}$ of closed $(1,1)$-forms on $M$. Locally, over $U_{\alpha}$ we have a spectral sheaf ${\mathcal S}_{\alpha}$, but now, as in Proposition \ref{sheafact}, tensoring with the line bundle $L_{\alpha\beta}$ gives an isomorphism from  ${\mathcal S}_{\alpha}$ to  ${\mathcal S}_{\beta}$ over $U_{\alpha}\cap U_{\beta}$. We deduce 
\begin{prp} A twisted generalized holomorphic bundle defines a spectral sheaf over the gerbe ${\mathcal G}$ on $TM$.
\end{prp}

The generic case consists of a spectral variety $S\subset TM$ and a trivialization of the gerbe ${\mathcal G}$. 

\begin{rmk} This holomorphic description in terms of the spectral variety is perhaps more transparent than  its direct image on $M$ -- from \cite{Bry} one can define  a direct image of a gerbe over ${\mathcal O}^*$ when the higher direct images $R^1\pi_*{\mathcal O}^*,R^2\pi_*{\mathcal O}^*$ vanish. In this case it is a gerbe over the sheaf of abelian groups $\pi_*{\mathcal O}^*$. 
\end{rmk}

\section{$\barD$-cohomology}  
\subsection{Hypercohomology groups}
In Section \ref{gen} we introduced generalized holomorphic bundles via an elliptic  complex:
$$\cdots \rightarrow C^{\infty}(V\otimes\Lambda^p E^{0,1})\stackrel {\barD}\rightarrow  C^{\infty}(V\otimes \Lambda^{p+1} E^{0,1})\rightarrow \cdots$$
We shall discuss here how the cohomology groups $H^q_{\bar D}(M,(V,\phi))$ of this complex can be determined in holomorphic terms. In particular, the interpretation above  of the twisted version as a sheaf on a gerbe presents a challenge, since such a sheaf has no obvious cohomology groups.

First consider an ordinary complex structure, where $E^{0,1}=T\oplus \barT^*$. Then 
$$\Lambda^r E^{0,1}=\bigoplus_{p+q=r}\Lambda^pT\otimes \Lambda^q\barT^*$$
The $\barD$-complex for a co-Higgs bundle $(V,\phi)$ is then  the sequence of operators 
$\bar\partial_A+(-1)^p\phi$ where
$$\bar\partial_A: \Omega^{0,q}(V\otimes \Lambda^pT)\rightarrow \Omega^{0,q+1}(V\otimes \Lambda^pT)$$
and
$$\phi: \Omega^{0,q}(V\otimes \Lambda^pT)\rightarrow \Omega^{0,q}(V\otimes \Lambda^{p+1}T)$$
acts as 
$$\phi^i\frac{\partial}{\partial z_i}(v\otimes W)=\phi^i(v)\frac{\partial}{\partial z_i}\wedge W$$

\begin{rmk} The index of the $\barD$-complex is the same as for $\phi=0$ by deformation invariance and this, by Riemann-Roch, is 
$$\sum_{p,q} (-1)^{p+q}\dim H^q(M,V\otimes \Lambda^pT)=\sum_{p=0}^m(-1)^p\ch(V)\ch(\Lambda^pT)\td(T)[M].$$
But $\Lambda^pT\cong \Lambda^{m-p}T^*\otimes K^*$ and from \cite{Hirz} (p. 92) 
$\sum (-1)^{p}\ch(\Lambda^p T^*)=(\td(T))^{-1}c_m(T)$
so the index is
\begin{equation}
\sum_{r=0}^{2m}(-1)^r\dim H_{\barD}^r(M,(V,\phi))=(-1)^m\ch(V)\ch(K^*)c_m(T)[M]=\rk V\chi(M)
\label{index}
\end{equation}
where $\chi(M)$ is the Euler characteristic. 
\end{rmk}

The operator $\bar\partial_A\pm \phi$ is the total differential in a double complex. If we pass from Dolbeault to \v Cech  cohomology it becomes the hypercohomology of the complex of sheaves
$$\cdots\rightarrow \mathcal O(V\otimes \Lambda^pT)\stackrel{\phi}\rightarrow \mathcal O(V\otimes \Lambda^{p+1}T)\rightarrow\cdots$$
so this is our $\barD$-cohomology.

First consider  the rank one case -- a line bundle $L$ and holomorphic vector field $X$. The complex of sheaves is 
$$\cdots\rightarrow \mathcal O_M(L\otimes \Lambda^pT)\stackrel{X\wedge}\rightarrow \mathcal O_M(L\otimes \Lambda^{p+1}T)\rightarrow\cdots$$
For a holomorphic vector bundle $W$, if $w$ is the tautological section of  $\pi^*W$ on $W$ then  we have the Koszul resolution of ${\mathcal O}_Z$ where $Z\subset W$ is the zero section:
$$\cdots\rightarrow \mathcal O_W( \Lambda^p \pi^*W^*)\stackrel{i_w}\rightarrow \mathcal O_W( \Lambda^{p-1}\pi^*W^*)\rightarrow\cdots$$
Take $W=T$, then identifying  $\Lambda^pT\cong \Lambda^{m-p}T^*\otimes K^*$, the exterior product by $X$ becomes the interior product. Our complex then becomes  the Koszul complex on $TM$ tensored with ${\mathcal O}_{X(M)}(LK^*)$ where $X:M\rightarrow TM$. It follows that $H^q_{\barD}(M,(L,X))$ is isomorphic (see for example \cite{GH}) to   $\Ext^q(X(M);{\mathcal O}_Z, LK^*).$ 

If the vector field $X$ has isolated zeros, then the section $X(M)$ and the zero section $Z$ intersect transversally  in a finite number of points and then from global duality  one can see that $\Ext^q(X(M);{\mathcal O}_Z, LK^*)$ vanishes if $q\ne m$, and  
 when $q=m$ is isomorphic to a space of sections of a line bundle over $X(M)\cap Z$. (Since $X$ vanishes at $\chi(M)=\#X(M)\cap Z$ points this is of course consistent with the index formula (\ref{index})).
 
 Since the twisted case  is not quite as standard,  we  give an elementary approach to this result for $X$ with isolated zeros  in terms of a Poincar\'e lemma for part of the $\barD$-complex. Recall that  
 $$\Lambda^r E^{0,1}=\bigoplus_{p+q=r}\Lambda^pT\otimes \Lambda^q\barT^*$$
 and denote by ${\mathcal E}^{r}$ the sheaf of $C^{\infty}$ sections of $ L\otimes\Lambda^r E^{0,1}$.
\begin{lemma}
  For $r>m$, and $r<m$,
  $0\rightarrow \barD{\mathcal E}^{r-1}\rightarrow {\mathcal E}^{r}\rightarrow  \barD{\mathcal E}^{r}\rightarrow 0$
  is exact. 
\end{lemma}
\begin{lemprf}

\noindent (i) Suppose that $r=m+k>m$ and $\varphi$ is a local section of  ${\mathcal E}^{r}$.  Then $\varphi= \varphi_{k}+\varphi_{k+1}+\dots+\varphi_m$ where $\varphi_p\in \Omega^{0,r-p}(L\otimes\Lambda^pT)$. Consider $\varphi_k\in  \Omega^{0,m}(L\otimes \Lambda^k T)$. Since $m$ is the top degree $\bar\partial \varphi_k=0$  trivially, and so $\varphi_k=\bar\partial\psi_k$ for some $\psi_k\in  \Omega^{0,m-1}(\Lambda^k T)$ by the Dolbeault lemma. Now $\barD \varphi=0$ implies that $\bar\partial\varphi_{k+1}\pm X\wedge\varphi_{k}=0$ so that  $\bar\partial(\varphi_{k+1}\pm X\wedge \psi_k)=0$ and hence 
$\varphi_{k+1}=\bar\partial \psi_{k+1}\mp X\wedge \psi_k.$
Continue 
until $\varphi_m\in \Omega^{0,k}(L\otimes \Lambda^mT)$, where since $k>0$ the lemma can still be applied.

\noindent (ii) Suppose $0<r<m$  then $\varphi=\varphi_0+\varphi_1+\dots+ \varphi_r$ and $\varphi_0\in \Omega^{0,r}(L)$. 
If $\barD\varphi=0$, then we have $\bar\partial \varphi_0=0$ and so if $r>0$, $\varphi_0=\bar\partial \psi_0$ where $\psi_0\in  \Omega^{0,r-1}(L)$. But $\barD\varphi=0$ implies that 
$\bar\partial\varphi_1+X\wedge \varphi_0=0$ and so 
$\bar\partial(\varphi_1+X\wedge \psi_0)=0$ and 
$\varphi_1=\bar\partial \psi_1-X\wedge \psi_0.$
Continue until the last term $\varphi_r$.

Now $0=\bar\partial \varphi_r\pm X\wedge\varphi_{r-1}=\bar\partial( \varphi_r\pm X\wedge \psi_{r-1})$ and so 
$\varphi_r\pm X\wedge \psi_{r-1}\in \Omega^{0,0}(L\otimes \Lambda^rT)$ is  a holomorphic section.  But $\barD\varphi=0$ implies also $X\wedge \varphi_r=0$ and so since $r<m$ by the Koszul resolution $\varphi_r\pm X\wedge \psi_{r-1}=X\wedge\rho_{r-1}$ for a holomorphic $\rho_{r-1}$, so $\varphi_r= X\wedge(\rho_{r-1}\mp\psi_{r-1})$. Putting all  terms together means that
$\varphi=\barD\psi$. 

\noindent (iii) When $r=0$, $\bar\partial\varphi_0=0$ so $\varphi_0$ is holomorphic. But 
$X\wedge \varphi_0=0$ and so $\varphi_0=0$ since $X$ is not identically zero.
\end{lemprf}

We now have: 
\begin{prp} Let $(L,X)$ be a generalized holomorphic line bundle over a compact $m$-dimensional complex manifold, where the vector field $X$ has non-degenerate zero set.  Then the $\barD$-cohomology vanishes except in degree $m$, and 
$H^m_{\barD}(M,(L,X))\cong H^0(X^{-1}(0),LK^*)$.
\end{prp}

\begin{prf} For $r<m$ the lemma and the usual \v Cech argument shows that the $\barD$-cohomology is isomorphic to the \v Cech cohomology of the sheaf defined by the kernel of $\barD:{\mathcal E}^0\rightarrow  {\mathcal E}^1$, but as in the lemma this sheaf is zero. In the case $r=m$, in the proof of the lemma  we have $\varphi_m\pm X\wedge \psi_{m-1}$ a holomorphic section of $L\otimes \Lambda^mT=LK^*$. We can solve for $\rho_{m-1}$ if and only if it  vanishes on the zero set of $X$.   For $r=m+k$   we find similarly at some stage a holomorphic $k$-cocycle with values in $LK^*$, but this can be changed by a coboundary to vanish on the zero set. So only the degree $m$ cohomology is nonvanishing and this is isomorphic to the sections of $LK^*$ on the zero set of $X$.
\end{prf}

In the general case we have the following:
\begin{thm} Let $(V,\phi)$ be a generalized holomorphic vector bundle over a compact $m$-dimensional complex manifold. Suppose the spectral sheaf is a line bundle $L$ supported on  a smooth $m$-dimensional submanifold $S\subset TM$ which  intersects the zero section $Z$ of $TM$ transversally and such that the projection to $\pi:S\rightarrow M$ is a finite map. 

Then the $\barD$-cohomology vanishes except in degree $m$, and 
$H^m_{\barD}(M,(V,\phi))\cong H^0(S\cap Z,L^*K_M^*)$. 
\end{thm}

\begin{rmk} Note here that $\pi$ expresses $S$ as a $k$-fold covering of $M$ where $\rk V=k$. Thus the homology class of $S$ is $k[M]$ and $\#(S\cap Z)=k[M].[M]=k\chi(M)$. This is in accordance with the vanishing and the index formula (\ref{index}).
\end{rmk}

\begin{prf} The vector bundle $V$ is the direct image of a line bundle $L'=L^*K_SK_M^*$ on $S$ and the Higgs field $\phi$  is the direct image of the tautological section $t$ of $\pi^*T$. Since $\pi$ is a finite map, the higher direct images of the sheaves on $S$ vanish. This means that 
$$\cdots\rightarrow \mathcal O(V\otimes \Lambda^pT)\stackrel{\phi}\rightarrow \mathcal O(V\otimes \Lambda^{p+1}T)\rightarrow\cdots$$
is a resolution. Repeating the argument above, we have $H^m_{\barD}(M,(V,\phi))\cong H^0(S\cap Z,L'K_M^*)$. But since $S$ and $Z$ are transversal we have a natural isomorphism on $S\cap Z$: $K_SK_M\cong K_{TM}\cong K_M^2$, so $K_S\cong K_M$  and $L'\cong L^*$. Hence 
$$H^m_{\barD}(M,(V,\phi))\cong H^0(S\cap Z,L^*K_M^*).$$
\end{prf}
\subsection{The twisted case}
In the twisted case there is no double complex to evaluate the cohomology. We can look at this from either the Dolbeault or \v Cech points of view. For the first, we can split the extension 
 $$0\rightarrow T^*\rightarrow E\stackrel{\pi}\rightarrow T\rightarrow 0$$
 to get an isomorphism $E\cong T\oplus T^*$. Then the Courant bracket on $E$ becomes twisted by a closed three-form $H$: $$[X+\xi,Y+\eta]_H=[X+\xi,Y+\eta]+i_Xi_YH.$$
 The splitting means that we can globally write 
 $\barD =\bar\partial_A\pm \phi$ where $\bar\partial_A:\Omega^{0,q}(V)\rightarrow \Omega^{0,q+1}(V)$. Now, however 
 $$\bar\partial_A^2=\pm i_{\phi}H^{1,2}$$
where $H^{1,2}$ is the $(1,2)$-component of $H$, so that $i_{\phi}H^{1,2}\in \Omega^{0,2}(\End V)$.

Alternatively, if we take a \v Cech point of view, we can argue that locally the picture is the standard co-Higgs case, and we have the Poincar\'e-type lemmas, but the identification is not consistent over threefold intersections. In fact the \v Cech differential $\delta$ satisfies
$$\delta^2=g$$
where $g=g_{\alpha\beta\gamma}$ is the holomorphic 2-cocyle defining the gerbe ${\mathcal G}$ (on $M$ in the line bundle case or $TM$ in general). 

In either case we have a total differential $\barD$ with $\barD^2=0$ but one of the constituent differentials has ``curvature".  However, recall from Section \ref{twistb} that the gerbe ${\mathcal G}$ was trivial on $Z$, so the spectral sheaf on the gerbe defines an ordinary sheaf on $Z$. More concretely, we write $B_{\alpha\beta}=\bar\partial \theta_{\alpha\beta}$ and then, using the tautological section $t$ of $\pi^*T$ on $TM$, $h_{\alpha\beta\gamma}= i_t(\theta_{\alpha\beta}+ \theta_{\beta\gamma}+ \theta_{\gamma\alpha})$  is a holomorphic function and  $g_{\alpha\beta\gamma}=\exp h_{\alpha\beta\gamma}$  a 2-cocycle which  defines the gerbe. Clearly $h_{\alpha\beta\gamma}$ is identically one on $Z$ since $t=0$ there, and so the spectral sheaf gives a well-defined line bundle $L$ over $S\cap Z$. It follows that:

\begin{thm} Let $(V,\phi)$ be a twisted generalized holomorphic vector bundle over a compact $m$-dimensional complex manifold. Suppose the spectral sheaf is a line bundle ${\mathcal L}$ over the  gerbe ${\mathcal G}$ supported on  a smooth $m$-dimensional submanifold $S\subset TM$ which  intersects the zero section $Z$ of $TM$ transversally and such that the projection to $\pi:S\rightarrow M$ is a finite map. 

Then the $\barD$-cohomology vanishes except in degree $m$, and 
$H^m_{\barD}(M,(V,\phi))\cong H^0(S\cap Z,L^*K_M^*)$ where  $L$ is the line bundle on $Z$ defined by  ${\mathcal L}$ and the canonical trivialization of ${\mathcal G}$ on $Z$.
\end{thm}

\begin{rmk} The situation above is actually irrelevant for line bundles in the K\"ahler case, because of the theorem of Carrell and Lieberman \cite{CL} (see also \cite{GH}) which shows that a compact K\"ahler manifold with a holomorphic vector field with isolated singularities has cohomology only of type $(p,p)$. 
Thus a nontrivial twist, which involves a class in $H^{1,2}$, is impossible. 
\end{rmk}

\vskip 1cm
 Mathematical Institute, 24-29 St Giles, Oxford OX1 3LB, UK
 
 hitchin@maths.ox.ac.uk

 \end{document}